\documentclass[12pt]{elsarticle}
\usepackage{amsfonts}
\usepackage[fleqn]{amsmath}
\usepackage{mathptmx}
\usepackage{mathtools,float}
\usepackage{latexsym,amsthm,amssymb}
\usepackage{mathrsfs,setspace,enumerate}
\usepackage[margin=0.7in]{geometry}
\usepackage{graphicx}
\usepackage{subcaption}
\usepackage[colorlinks]{hyperref}
\usepackage{float}
\usepackage{multirow}
\usepackage{verbatim}
\hypersetup{
	colorlinks=true,
	linkcolor=blue,
	citecolor=green,
	filecolor=magenta,
	urlcolor=blue
}

\makeatletter
\def\ps@pprintTitle{%
	\let\@oddhead\@empty
	\let\@evenhead\@empty
	\def\@oddfoot{\reset@font\hfil\thepage\hfil}
	\let\@evenfoot\@oddfoot
}
\makeatother
\begin{document}
\begin{abstract}
	This manuscript focuses on the construction of compactly supported dual Gabor frames in $L^2(\mathbb{R})$. The performance of the constructed dual frames is analyzed for Gabor systems generated by B-splines and exponential B-splines of orders 2 and 3. The reconstruction performance of these dual windows is evaluated using the average mean square error (AMSE) for standard one-dimensional benchmark signals. For two-dimensional data, image reconstruction experiments are carried out using tensor-product Gabor frames, and the reconstruction accuracy is also assessed using AMSE.
		Using the duality condition for Gabor systems \cite{jan}, several alternate dual windows with finite support are constructed under suitable assumptions such as the partition of unity property. Additional dual windows can also be obtained from an existing dual. The canonical dual window admits an explicit expression that avoids direct inversion of the frame operator and yields reconstruction errors close to numerical precision. The constructed non-canonical compactly supported duals also exhibit stable and competitive reconstruction performance. These findings indicate that compactly supported dual windows based on B-spline and exponential B-spline generators provide effective and practical alternatives for signal and image processing applications, particularly in situations where compact support and computational efficiency are important.
\end{abstract}	\begin{keyword}
Gabor frame, B-splines, Exponential B-splines, Dual frames
\MSC[2020] 42C15 \sep 42C40. 
	\end{keyword}
	
	\begin{frontmatter}
	\title{Performance Assessment and Construction of Compactly Supported Dual Windows for B-spline and Exponential B-spline Gabor Frames}
		\author[1]{Sruthi Raghoothaman}
		\ead{sruthiragoothaman6@gmail.com}
		\author[2]{Noufal Asharaf}
		\ead{noufal@cusat.ac.in}
		\address[1]{Department of Mathematics, Cochin University of Science and Technology}
	\end{frontmatter}
	
	\section{Introduction}
Gabor frames play a fundamental role in time–frequency analysis and signal processing, as they provide stable and redundant representations of functions in $L^2(\mathbb{R})$. 
They arise from the action of time and frequency shift operators on a fixed window function.
  For $a,b >0$, and $g \in L^2(\mathbb{R})$ let $T_a$ and $E_b$ denote the translation (time-shift) and modulation (frequency-shift) operators on $L^2(\mathbb{R})$,  defined by $
T_a f(x)=f(x-a)$, and 
$E_b f(x)=e^{2\pi i b x} f(x),\  x\in\mathbb{R}.
$
For $g \in L^2(\mathbb{R})$, the collection of time-frequency shifts of $g$ given by
$\mathcal{G}(g,a,b) = \{E_{mb}T_{na}g\}_{m,n \in \mathbb{Z}}$
is called a Gabor system generated by $g$ with time shift $a$ and frequency shift $b$. The system $\mathcal{G}(g,a,b)$ is called a Gabor frame for $L^2(\mathbb{R})$ if there exist constants $A,B>0$ such that:
$$ A \|f\|^2 \leq \sum_{m,n \in \mathbb{Z}} \left| \langle f, E_{mb} T_{na} g \rangle \right|^2 \leq B \|f\|^2, \ \ \forall f \in L^2(\mathbb{R}).$$
$\mathcal{G}(g,a,b)$ is a Bessel sequence in $L^2(\mathbb{R})$
 if at least the upper inequality holds. It is well known that for the Gabor frame $\mathcal{G}(g,a,b)$, there exists a function 
 $h \in L^2(\mathbb{R})$ such that every $f \in L^2(\mathbb{R})$ admits the reconstruction 
 $
 f = \sum_{m,n\in\mathbb{Z}} 
 \langle f, E_{mb}T_{na}h \rangle \, E_{mb}T_{na}g.
 $
Due to the redundancy of frames, each function 
$f \in L^2(\mathbb{R})$ typically possesses infinitely many frame representations. 
In particular, when the lattice parameters satisfy $ab<1$, the associated Gabor system is overcomplete, and we can find infinitely many such functions h in $L^2(\mathbb{R})$. Any such function, with the
additional property that $ \mathcal{G}(g,a,b) $  satisfies the upper frame condition,
is called a dual generator. We refer $(g, h)$ as a pair of dual frame generators.
 Among all possible dual windows, a classical choice is 
 $h = S^{-1}g$, where $S$ denotes the frame operator associated with 
 $\mathcal{G}(g,a,b)$.
   The frame operator of $ \mathcal{G}(g,a,b) $ is given by $$ S f := \sum_{m,n \in \mathbb{Z}} \langle f, E_{mb} T_{na} g \rangle E_{mb} T_{na} g $$
    The frame operator is bounded, bijective, self-adjoint, and positive. The frame $ \{ E_{mb}T_{na}S^{-1}g\}_{m,n \in \mathbb{Z}}$ called the canonical dual. Several characterizations of Gabor frames have been developed in the literature, notably by Wexler and Raz \cite{wexler}, Daubechies et al. \cite{daub}, and Ron and Shen \cite{shen}.
	
The canonical dual of a Gabor frame plays a central role in exact signal reconstruction. 
It is obtained by applying the inverse of the associated frame operator to the generating window, and it guarantees perfect reconstruction of all functions in $L^2(\mathbb{R})$. 
However, the canonical dual does not necessarily inherit desirable properties of the original window, such as smoothness, compact support, or favorable decay in time or frequency. 
In particular, even when the generating window is compactly supported, its canonical dual typically has infinite support. This limitation motivates the search for alternative (non-canonical) dual windows possessing improved structural or practical properties. 
Since compact support is closely related to computational efficiency, numerical stability, and localized implementations, the construction of compactly supported dual Gabor windows constitutes an important problem in Gabor analysis. 
It is worth noting that any dual window depends on both the original window $g$ and the lattice parameters $a$ and $b$.
\par In this article, we develop constructive methods for obtaining compactly supported dual windows for selected Gabor generators. 
From these initial constructions, we derive further families of dual windows. 
Explicit examples are provided, and the performance of the proposed duals is evaluated in the context of signal and image reconstruction. In Section 2, we review several well-known constructions used for obtaining dual windows. In Section 3, we summarize a general approach for constructing alternative duals for a given frame, as proposed in \cite{arefi1}, along with the complete characterization of compactly supported dual windows given in \cite{Sdiana}. These construction techniques form the basis for the numerical experiments conducted in this study. Section 4 provides a brief overview of the Gabor window functions used in the experiments: the symmetric B-splines of orders 2 and 3, and the exponential B-splines of orders 2 and 3.For each window, several associated dual windows are constructed using the methods developed in the preceding sections, with the translation and modulation parameters $a$ and $b$ fixed throughout the experiments. A comparison of reconstruction performance are done based on Average Mean Squared Error (AMSE). Finally, Section 5 investigates the effectiveness of the constructed dual windows in practical image processing applications and provides a comparative analysis of their performance.
\section{Dual Gabor Frame Constructions with Compact Support}
Several duality principles in Gabor frame theory have been developed in the literature; see, for example, Janssen \cite{aug}, Grochenig \cite{groch}, Casazza et al. \cite{casazza}, Christensen \cite{olepairs}, and Christensen and Kim \cite{kim}. Among these results, the duality condition given by Janssen plays a key role in the construction of compactly supported Gabor duals, as it provides a necessary and sufficient characterization for two Gabor systems to form a dual pair. In particular, this condition translates the duality problem into a localized functional identity involving the generating windows, thereby enabling explicit constructions. More precisely, two Bessel sequences $\mathcal{G}(g,a,b)$ and $\mathcal{G}(h,a,b)$ form dual frames for $L^2(\mathbb{R})$ if and only if for $m\in \mathbb{Z}$ \cite{jan},	\begin{equation*}\label{duality}
	\sum_{n \in \mathbb{Z}} g(x -\frac{m}{b} + na)\,h(x + na) = b\delta_{m,0}, \quad \text{for a.e. } x \in [0, a].
\end{equation*} 
This condition represents an infinite system of equations, each 
involving infinitely many terms. Consequently, in order to explicitly 
determine a dual window $h$, it is necessary to impose suitable 
assumptions on the parameters $a$, $b$, and the generating function $g$.  Throughout the remainder of the paper, we assume that 
$g \in L^2(\mathbb{R})$ is compactly supported in 
$\left[-\frac{N}{2}, \frac{N}{2}\right]$ for some $N\in \mathbb{N}$. Furthermore, we restrict 
our attention to dual windows that can be written as finite linear 
combinations of integer translates of $g$. In our construction, we restrict attention to constructing dual windows that can be expressed as finite linear combinations of integer shifts of the generating window $g$.  
More precisely, we consider dual windows of the form
\begin{equation*}
h(x) = \sum_{n=-N+1}^{N-1} a_n \, g(x+n),
\end{equation*}
for suitable real coefficients $\{a_n\}_{n=-N+1}^{N-1}$. 
This structure allows us to preserve several desirable properties of the original window. 
In particular, compact support and regularity properties of $g$ are automatically inherited by $h$.
Moreover, if $g$ belongs to the Schwartz space $\mathcal{S}(\mathbb{R})$, the Wiener space, or Feichtinger’s algebra $S_0(\mathbb{R})$, then any dual window of the above form belongs to the same function space \cite{olebook}.
For simplicity, we restrict to the case $a=1$. 
Assume that $g \in L^2(\mathbb{R})$ is real-valued, bounded, and satisfies 
$
\text{supp}(g) \subset [0,N],
$
together with the partition of unity property
$
\sum_{n\in\mathbb{Z}} g(x-n) = 1, 
\quad x\in\mathbb{R}.
$

Let $b \in \left(0,\frac{1}{2N-1}\right]$, and consider real coefficients $\{a_n\}_{n=-N+1}^{N-1}$ satisfying
\[
a_0 = b, 
\qquad 
a_n + a_{-n} = 2b, 
\quad n=1,2,\dots,N-1.
\]
Then 
$\mathcal{G}(g,1,b)$ and $\mathcal{G}(h,1,b)$
generates a pair of dual Gabor frames for $L^2(\mathbb{R})$, see Theorem 3.1 of \cite{kim}. Theorem 3.1 is stated for functions $g$ with $\text{supp } g \subset [0,N]$, but it remains valid if 
$\text{supp}(g)\subset [-N/2,N/2]$.
 A particular choice of the coefficients yields explicit symmetric and asymmetric dual windows. 
If we set 
$
a_n = b,$ for $ n=-N+1,\dots,N-1,
$
then the resulting dual window takes the symmetric form
\begin{equation}\label{symdual}
	h(x) = b \sum_{n=-N+1}^{N-1} g(x+n).
\end{equation}
The function satisfies $h=b$ on the support of $g$, [see Corollary 3.2 in \cite{kim}].
Alternatively, choosing the coefficients so that
$
a_0=b, \ a_n=2b \text{ for } n=1,\dots,N-1, 
\  a_{-n}=0 \text{ for } n=1,\dots,N-1,
$
produces the asymmetric dual window [see Theorem 2.2 in \cite{olepairs},]
\begin{equation}\label{asymdual}
	k(x) = b g(x) + 2b \sum_{n=1}^{N-1} g(x+n),
\end{equation}
which has the shortest possible support among dual windows of this form.

The construction method described above applies to compactly supported 
windows supported on $[-\frac{N}{2},\frac{N}{2}]$ with the partition of unity, and for lattice parameters 
$a=1$ and $b\in \left(0,\frac{1}{2N-1}\right]$. 
Under these assumptions, the windows $g$ and $h$ do not in general generate 
dual Gabor frames when $a\neq 1$.
However, this restriction can be overcome by a suitable dilation. 
Let $D_a$ denote the dilation operator defined by 
$
D_a f(x) = a^{-1/2} f(x/a), \  a>0.
$
If $a,b>0$ satisfy $ab\in\left(0,\frac{1}{2N-1}\right]$, then the dilated windows 
$D_a g$ and $D_a h$ generate a pair of dual Gabor frames 
$
\mathcal{G}(D_a g,a,b) \text{ and } 
\mathcal{G}(D_a h,a,b)
$
for $L^2(\mathbb{R})$, where $
h(x) = ab\, g(x) + 2ab \sum_{n=1}^{N-1} g(x+n),
$ see Theorem 2.6 of \cite{olepairs}.

\section{Construction of Dual Windows from Known Duals}

Additional compactly supported dual windows can be obtained from a given dual window. 
According to Corollary 2.2 in \cite{arefi1}, if $g^d$ is a dual window of $g$, then alternate dual windows can be constructed via
\begin{equation*}\label{altdual}
	\tilde{g} = S^{-1}g - g + S g^d,
\end{equation*}
where $S$ denotes the frame operator associated with $\mathcal{G}(g,a,b)$. In particular, starting from an initial dual window $\tilde{g}_0 = g^d$, one can define recursively
\begin{equation}\label{altdual*}
	\tilde{g}_{l+1} = S^{-1}g - g + S \tilde{g}_l,
\end{equation}
and each function $\tilde{g}_l$ thus obtained is again a dual window of $g$.
As discussed earlier, computing the canonical dual window $S^{-1}g$ is, in general, a nontrivial task, since it requires inversion of the frame operator. 
However, in certain situations, an explicit formula for $S^{-1}g$ is available, which makes the above construction practically implementable. Since the window $g$ we consider here is compactly supported with support length $N$, and since the function
\[
G(x) = \sum_{n \in \mathbb{Z}} |g(x - na)|^2
\]
is bounded and bounded away from zero for $a=1$, it follows from Lemma 1.2 in \cite{olepairs} that, for $b \le \frac{1}{N}$, the Gabor system $\mathcal{G}(g,a,b)$ forms a frame for $L^2(\mathbb{R})$. Moreover, the associated canonical dual window is given explicitly as
\begin{equation}\label{cano}
	S^{-1}g(x) = \frac{b}{G(x)}\, g(x).
\end{equation}
Since $g$ has compact support, only finitely many terms contribute to the function $G$ at each point. Moreover, $S^{-1}g$ is supported on the same interval as $g$, hence the canonical dual is also compactly supported. 
If the frame is not tight  (i.e., $A\neq B$), then the factor $1/G(x)$ is not constant, so the expression for $S^{-1}g$ is not explicit in a constructive sense. 
While the formula tells us the value of $S^{-1}g$ pointwise, it does not immediately yield information about the decay of its Fourier transform. Information about the decay properties can be obtained via the smoothness of the function.

\par  A constructive method for finding all compactly supported dual Gabor frames is described in \cite{Sdiana}. The method offers several significant advantages over traditional approaches. Unlike general existence results, this method provides explicit formulas for dual windows under clear and verifiable conditions, making it practical for implementation. A significant advantage of the method is its applicability to compactly supported generators, which enables the construction of dual windows that retain compact support—a feature crucial in applications requiring time localization.
Suppose that $g,g^d\in L^2(\mathbb{R})$ and $a,b>0$ are such that 
$\mathcal{G}(g,a,b)$ is a Gabor frame for $L^2(\mathbb{R})$ and 
$\mathcal{G}(g^d,a,b)$ is a dual frame of $\mathcal{G}(g,a,b)$. 
Then every dual window of $\mathcal{G}(g,a,b)$ can be obtained from $g^d$ by adding a correction term determined by 
\begin{equation*}\label{stoeva}
	h	= g^d + w - \sum_{m,n\in\mathbb{Z}}
	\langle g^d, E_{mb}T_{na}g\rangle
	\, E_{mb}T_{na}w,
\end{equation*}
where $u\in L^2(\mathbb{R})$ is chosen such that 
$\mathcal{G}(w,a,b)$ forms a Bessel sequence in $L^2(\mathbb{R})$.
 We observe that if $g$ and $g^d$ are compactly supported, then there exist a finite set $K$ (depends on $g$, $g^d$ and $a$) such that $\langle g^d, E_{mb}T_{na}g \rangle \neq 0 $ if and only if $n\in K.$
For implementation purposes, it is convenient to use real-valued windows. 
A more explicit characterization can be obtained in the real-valued setting. 
Assume in addition that $g$ and $g^d$ are real-valued. 
For $m,n\in\mathbb{Z}$, define the operators
\[
P_{m,n}(f)(x)
:=
\left(
\int g^d(x)\, f(2\pi m b x)\, T_{na}g(x)\, dx
\right)
f(2\pi m b x).
\]
Then every real-valued compactly supported dual window of $\mathcal{G}(g,a,b)$ can be written in the form
\begin{equation}\label{stoevaReal}
\varphi
=
g^d
+
w
-
\sum_{n\in K}
\left(
\langle g^d, T_{na}g\rangle
+
2\sum_{m\in\mathbb{N}}
\big(
P_{m,n}(\cos)
+
P_{m,n}(\sin)
\big)
\right)
T_{na}w,
\end{equation}
where $w\in L^2(\mathbb{R})$ is real-valued, compactly supported, and chosen such that the Gabor system 
  $\mathcal{G}(w,a,b)$ 
forms a Bessel sequence in $L^2(\mathbb{R})$.
 
\par
 We provide a detailed comparison of these two dual construction methods for different Gabor windows in the next section. Both \cite{Sdiana} and \cite{arefi1} address the problem of constructing alternate dual windows for Gabor frames, but they approach the problem from distinct perspectives. In \cite{arefi1}, the focus is on an explicit algebraic method for constructing alternate duals by perturbing the canonical dual through linear combinations of shifts of the original window. Their iterative formula, equation \eqref{altdual*}, can be used in practical situations only when an explicit expression for the inverse of the frame operator is available, as in \eqref{cano}. In contrast, \cite{Sdiana} presents a more general and structural characterization of all compactly supported dual windows of a Gabor frame without relying on operator inversion. This approach is versatile, applying to a broad class of Gabor frames, and is particularly suited for scenarios where the canonical dual is not explicitly known or lacks desired localization properties. 
\section{Dual Construction for B-Splines and Exponential B-Splines} 
\par For the numerical computations, we use B-splines of orders 2 and 3, denoted by $B_2$ and $B_3$, and exponential B-splines of orders 2 and 3, denoted by $\varepsilon_2$ and $\varepsilon_3$ as Gabor windows. These functions are widely used in time-frequency analysis due to their favorable localization properties, smoothness, and ease of implementation. The B-splines $B_2$ and $B_3$
are piecewise polynomial functions with compact support, making them computationally efficient while still providing good overlap properties required for frame constructions. Unlike polynomial B-splines, exponential B-splines are defined with an exponential weight that allows them to adapt more effectively to exponentially decaying signals. This makes them particularly useful in applications where the underlying data exhibits non-polynomial behavior. In the context of Gabor frames, exponential B-splines can serve as effective generating windows because they provide smoothness, compact support, and good localization in both time and frequency. These properties ensure stable and accurate reconstruction of signals from frame coefficients.

Symmetric B-splines are defined inductively by
\[
B_1 = \chi_{\left[-\frac{1}{2},\,\frac{1}{2}\right]},
\qquad
B_{N+1} = B_N * B_1,
\]
where $*$ denotes convolution. 
It follows that $B_N$ is supported on the interval 
$
\left[-\frac{N}{2},\,\frac{N}{2}\right].
$
Moreover, the integer translates of $B_N$ satisfy the partition of unity property,
$
\displaystyle \sum_{n\in\mathbb{Z}} B_N(x-n) = 1,
\  x\in\mathbb{R}.
$ We now recall the explicit formulas for the second- and third-order symmetric B-splines.

The second-order symmetric B-spline $B_2$ is given by
\[
B_2(x)=
\begin{cases}
	1+x, & x\in[-1,0], \\[1ex]
	1-x, & x\in[0,1], \\[1ex]
	0, & \text{otherwise}.
\end{cases}
\]

\begin{figure}[H]\label{B2}
	\centering
	
	\begin{subfigure}[b]{0.3\textwidth}
		\centering
		\includegraphics[width=\linewidth]{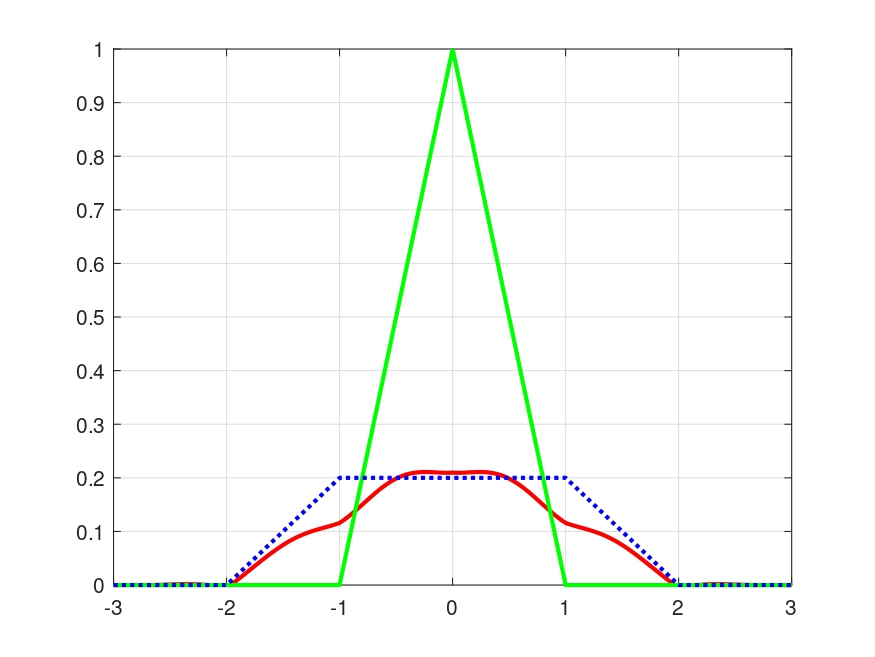}
		\subcaption{}
		\label{fig:sub-a}
	\end{subfigure}
	\hfill
	\begin{subfigure}[b]{0.3\textwidth}
		\centering
		\includegraphics[width=\linewidth]{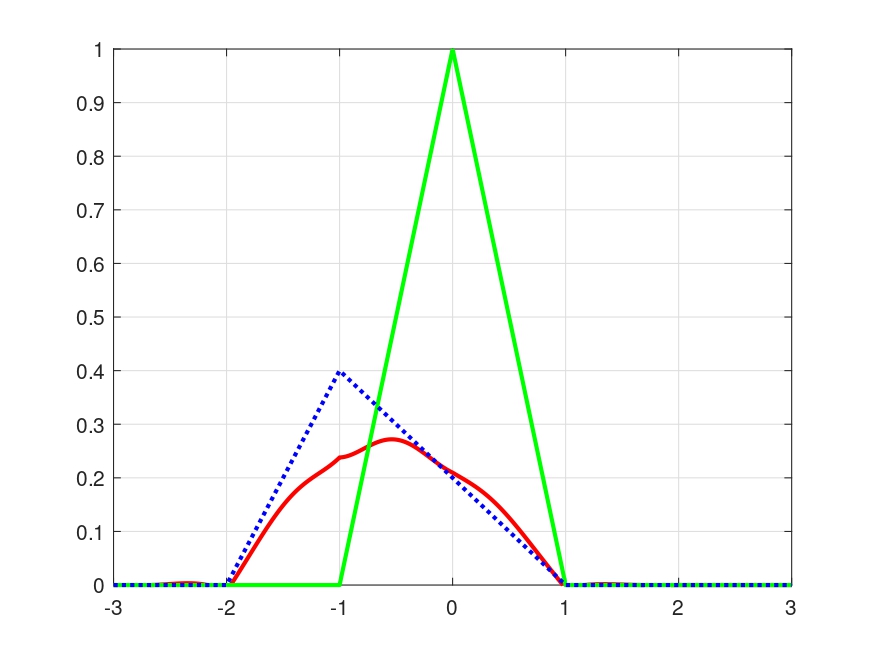}
		\subcaption{}
		\label{fig:sub-b}
	\end{subfigure}
	\hfill
	\begin{subfigure}[b]{0.3\textwidth}
		\centering
		\includegraphics[width=\linewidth]{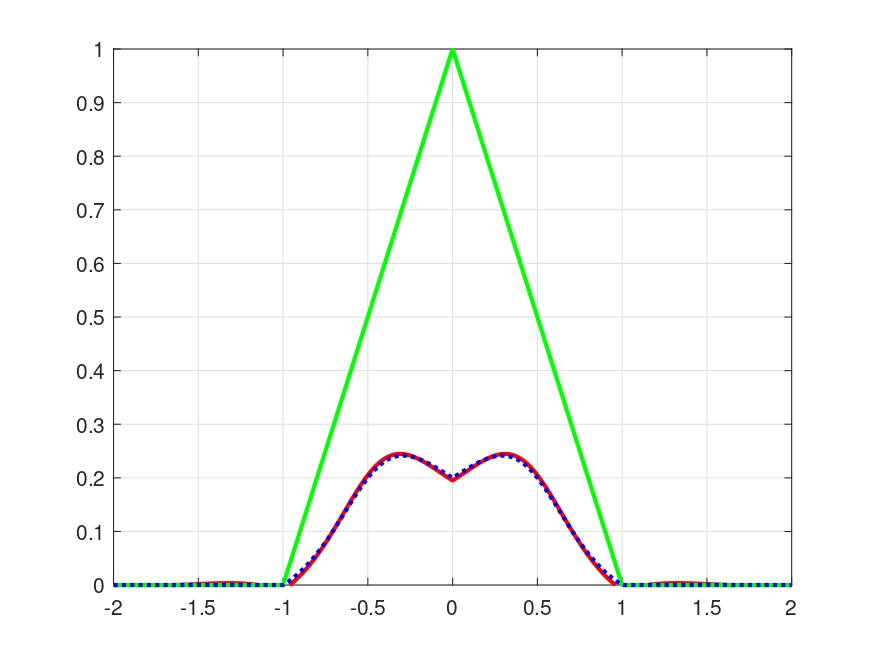}
		\subcaption{}
		\label{fig:sub-c}
	\end{subfigure}
	
	\caption{$a = 1$, $b = \frac{1}{5}$. Green curve represents the generator $B_2$. 
		(a) Blue: symmetric dual $k$; red: the dual constructed from $k$ using \ref{stoevaReal}.
		(b) Blue: asymmetric dual $h$; red: the dual constructed from $h$ using using \ref{stoevaReal}. 
		(c) Blue: canonical dual $S^{-1}B_2$; red: the dual constructed from it using using \ref{stoevaReal}.}
	\label{fig:B2}
\end{figure}

\begin{figure}[H]\label{B3}
	\centering
	
	\begin{subfigure}[b]{0.3\textwidth}
		\centering
		\includegraphics[width=\linewidth]{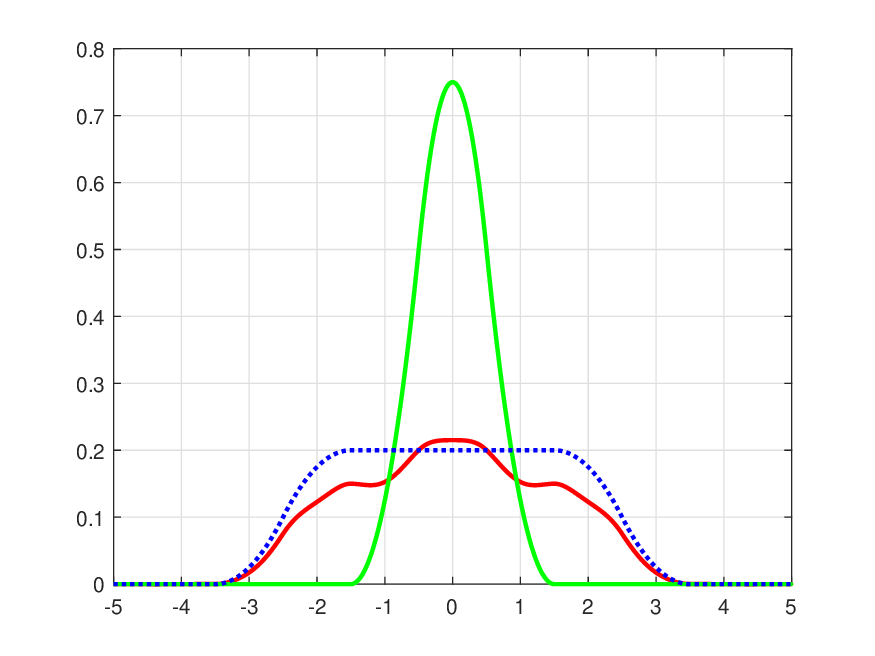}
		\subcaption{}
		\label{fig:sub-a}
	\end{subfigure}
	\hfill
	\begin{subfigure}[b]{0.3\textwidth}
		\centering
		\includegraphics[width=\linewidth]{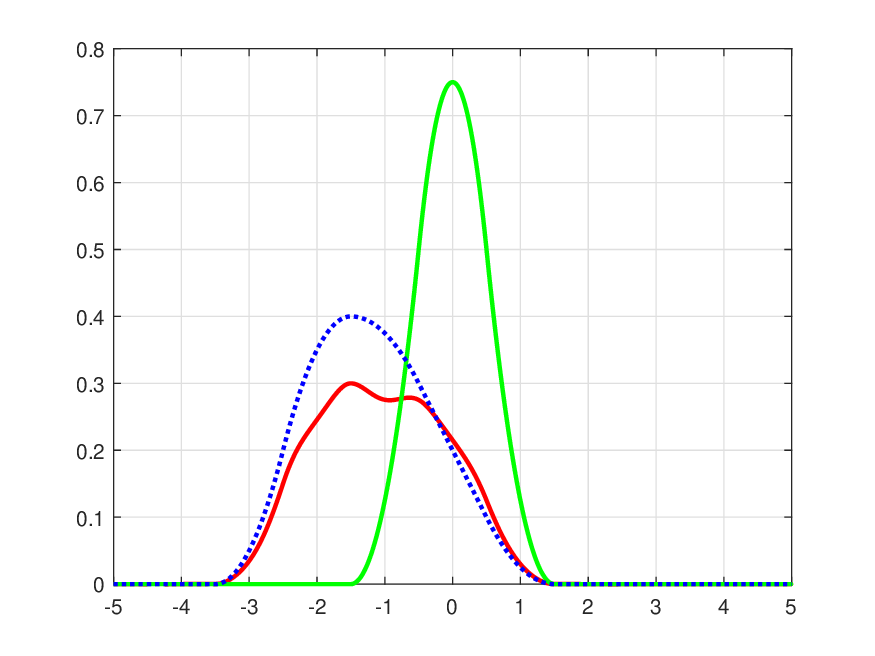}
		\subcaption{}
		\label{fig:sub-b}
	\end{subfigure}
	\hfill
	\begin{subfigure}[b]{0.3\textwidth}
		\centering
		\includegraphics[width=\linewidth]{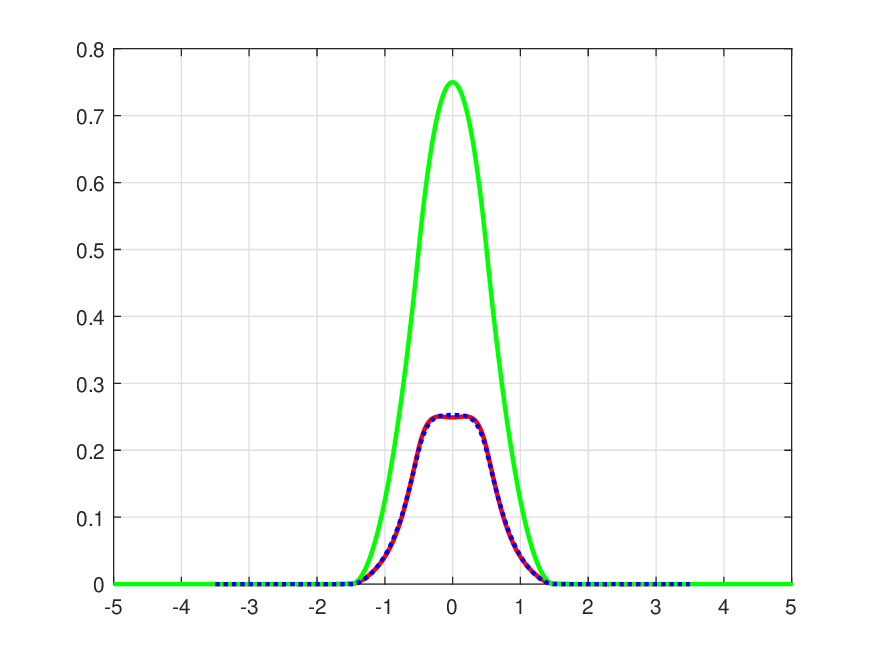}
		\subcaption{}
		\label{fig:sub-c}
	\end{subfigure}
	
	\caption{$a = 1$, $b = \frac{1}{5}$. Green curve represents the generator $B_3$. 
		(a) Blue: symmetric dual $k$; red: the dual constructed from $k$ using \ref{stoevaReal}. 
		(b) Blue: asymmetric dual $h$; red: the dual constructed from $h$ using  \ref{stoevaReal}. 
		(c) Blue: canonical dual $S^{-1}B_3$; red: the dual constructed from it using \ref{stoevaReal}.}
	\label{fig:B3}
\end{figure}

The third-order symmetric B-spline $B_3$ has the piecewise quadratic representation
\[
B_3(x)=
\begin{cases}
	\displaystyle \frac{1}{2}x^2 + \frac{3}{2}x + \frac{9}{8}, 
	& x\in\left[-\frac{3}{2},-\frac{1}{2}\right], \\[2ex]
	
	\displaystyle -x^2 + \frac{3}{4}, 
	& x\in\left[-\frac{1}{2},\frac{1}{2}\right], \\[2ex]
	
	\displaystyle \frac{1}{2}x^2 - \frac{3}{2}x + \frac{9}{8}, 
	& x\in\left[\frac{1}{2},\frac{3}{2}\right], \\[2ex]
	
	0, & \text{otherwise}.
\end{cases}
\]
Figures \ref{fig:B2} and \ref{fig:B3} show $B_2$ and $B_3$ along with their duals obtained from \eqref{stoevaReal}.

Let \( n \in \mathbb{N} \) and let \( a := (a_1, \ldots, a_N) \), where \( a_1, \ldots, a_N \in \mathbb{R} \) with 
\( a_i \neq 0 \) for at least one \( i \in \{1, \ldots, N\} \).
 An exponential B-spline 
\( \varepsilon_{N,a}' : \mathbb{R} \to \mathbb{R} \) of order \( N \) for the \( N \)-tuple \( a \) is a function of the form
\[
\varepsilon_{N,a}' := e^{a_1(\cdot)}\chi * e^{a_2(\cdot)}\chi * \cdots * e^{a_N(\cdot)}\chi,
\]
For $a_1 < a_2 < \cdots < a_N$, with $N \geq 2$, an explicit formula for $\varepsilon_N'$ is given by [Theorem 2.2,\cite{peter}],
\[
\varepsilon_N'(x) = 
\begin{cases}
	0, & x \notin [0,N], \\[1.2ex]
	
	\displaystyle 
	\sum_{i=1}^N 
	\left( \prod_{\substack{j=1 \\ j\neq i}}^N 
	\frac{1}{a_i-a_j} \right) 
	e^{a_i x},
	& x \in [0,1], \\[3ex]
	
	\displaystyle 
	(-1)^{k-1}
	\sum_{i=1}^N 
	\Bigg(
	\sum_{\substack{
			1\le j_1<\cdots<j_{k-1}\le N\\
			j_1,\dots,j_{k-1}\neq i}}
	\!\!\!
	\left(
	e^{a_{j_1}} + \cdots + e^{a_{j_{k-1}}}
	\right)
	\prod_{\substack{j=1 \\ j\neq i}}^N 
	\frac{1}{a_i-a_j}
	\Bigg)
	e^{a_i(x-k+1)},
	& x \in [k-1,k], \\[1ex]
	& k=2,\ldots,N .
\end{cases}
\]
Unlike classical B-splines, exponential B-splines do not, in general, satisfy the partition of unity property. In particular, it is not always true that
$
\displaystyle \sum_{n\in\mathbb{Z}} \varepsilon_N'(x-n) = 1, \  x \in \mathbb{R}.
$
However, certain exponential B-splines satisfy a weaker property of the form
$
\displaystyle \sum_{n\in\mathbb{Z}} \varepsilon_N'(x-n) = C, \  x \in \mathbb{R},
$
where \(C\) is a constant independent of \(x\).
Moreover, if all exponential parameters satisfy \(a_i \neq 0\) for every \(i\), then the partition of unity property cannot be achieved. By Theorem 3.1 of \cite{peter}, $\varepsilon_2'$ and $\varepsilon_3'$ will have this property for proper choice of exponential parameters. 
For $a_1 = p\neq 0$ and $a_2=0$,
\begin{equation*}\label{six}
	\varepsilon_2'(x) =
	\begin{cases}
		0, & x \notin [-1,1], \\[1ex]
		\displaystyle \frac{e^{p\left(x+\frac{1}{2}\right)} - e^{-\frac{p}{2}}}{p}, 
		& x \in [-1,0], \\[2ex]
		\displaystyle \frac{e^{\frac{p}{2}} - e^{p\left(x-\frac{1}{2}\right)}}{p}, 
		& x \in [0,1].
	\end{cases}
\end{equation*}

\begin{figure}[H]\label{EBS2}
	\centering
	
	\begin{subfigure}[b]{0.3\textwidth}
		\centering
		\includegraphics[width=\linewidth]{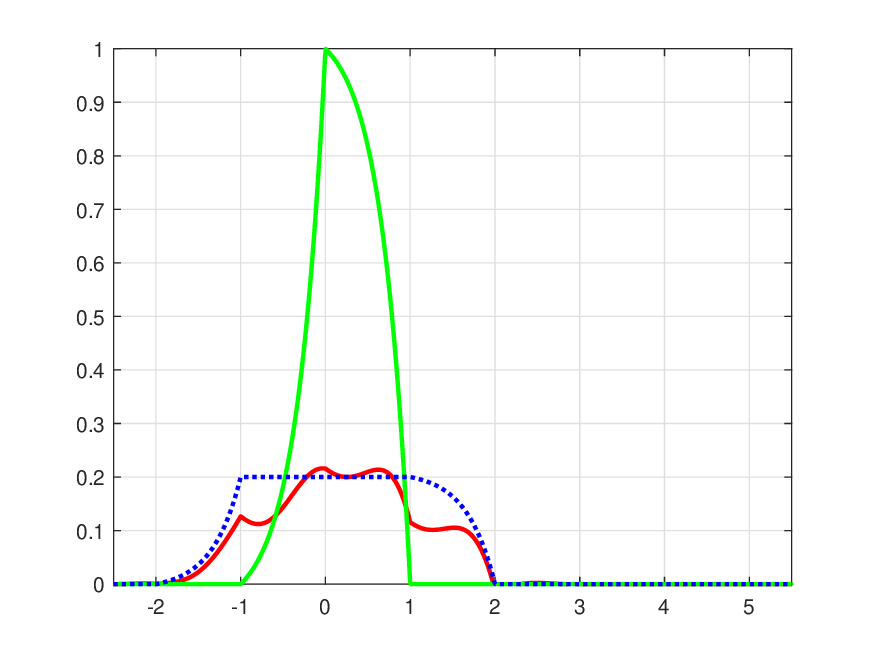}
		\subcaption{}
		\label{fig:sub-a}
	\end{subfigure}
	\hfill
	\begin{subfigure}[b]{0.3\textwidth}
		\centering
		\includegraphics[width=\linewidth]{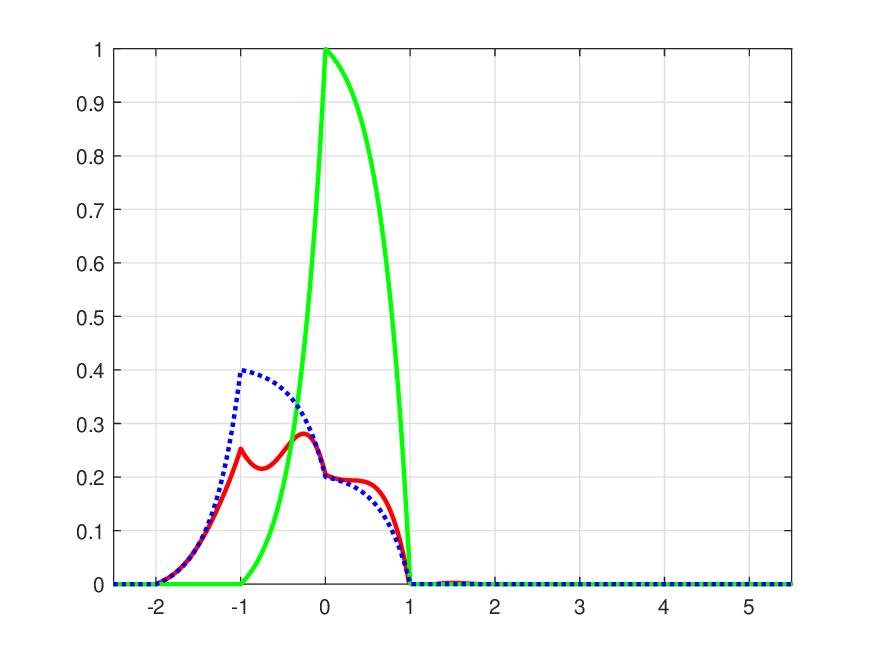}
		\subcaption{}
		\label{fig:sub-b}
	\end{subfigure}
	\hfill
	\begin{subfigure}[b]{0.3\textwidth}
		\centering
		\includegraphics[width=\linewidth]{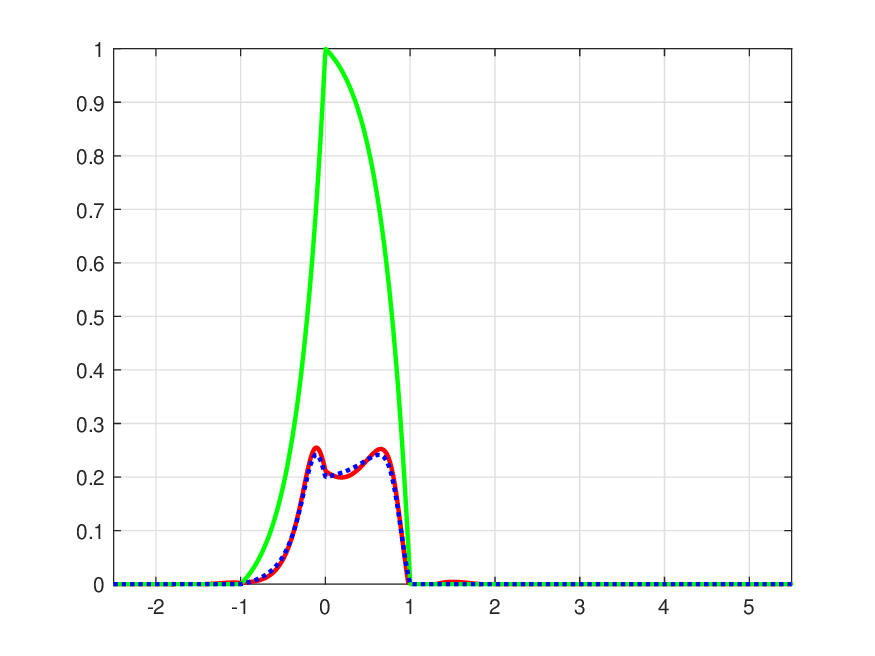}
		\subcaption{}
		\label{fig:sub-c}
	\end{subfigure}
	
	\caption{$a = 1$, $b = \frac{1}{5}$. Green curve represents the generator $\varepsilon_3$. 
		(a) Blue:  dual $k$; red: the dual constructed from $k$ using \ref{stoevaReal}.
		(b) Blue:  dual $h$; red: the dual constructed from $h$ using \ref{stoevaReal}. 
		(c) Blue: canonical dual $S^{-1}\varepsilon_3$; red: the dual constructed from it using \ref{stoevaReal}.}
	\label{fig:EBS2}
\end{figure}

\begin{figure}[H]\label{EBS3}
	\centering
	
	\begin{subfigure}[b]{0.3\textwidth}
		\centering
		\includegraphics[width=\linewidth]{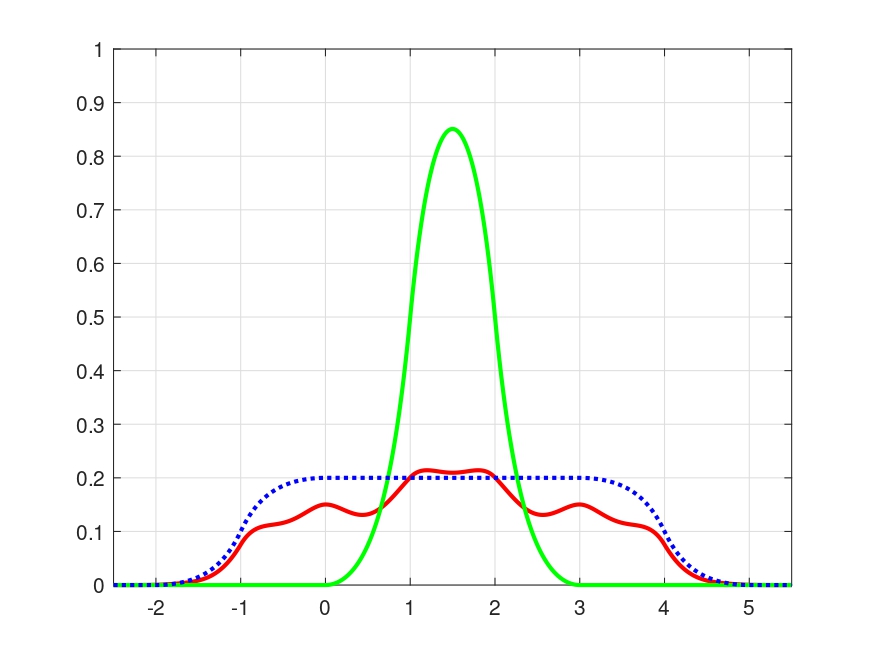}
		\subcaption{}
		\label{fig:sub-a}
	\end{subfigure}
	\hfill
	\begin{subfigure}[b]{0.3\textwidth}
		\centering
		\includegraphics[width=\linewidth]{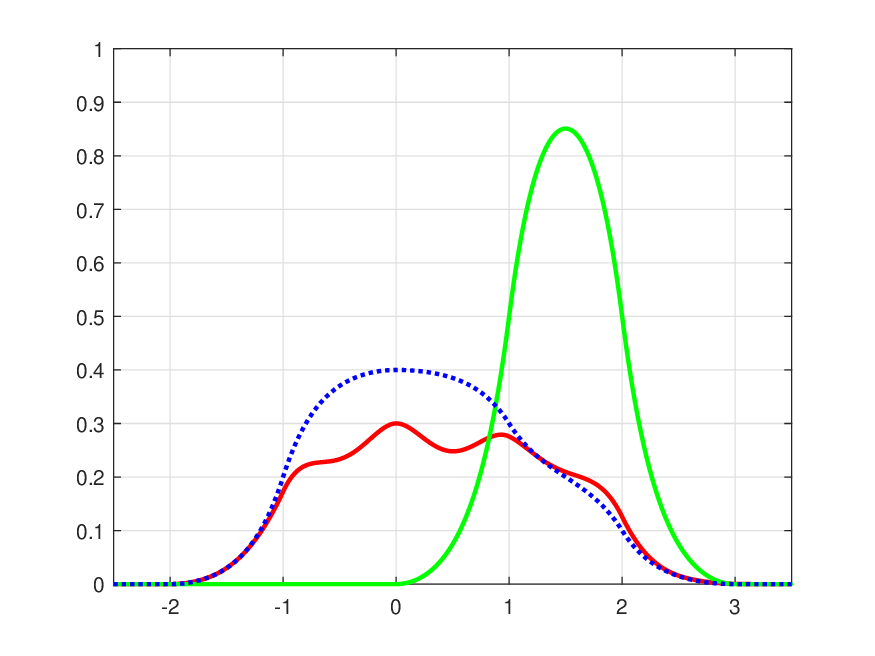}
		\subcaption{}
		\label{fig:sub-b}
	\end{subfigure}
	\hfill
	\begin{subfigure}[b]{0.3\textwidth}
		\centering
		\includegraphics[width=\linewidth]{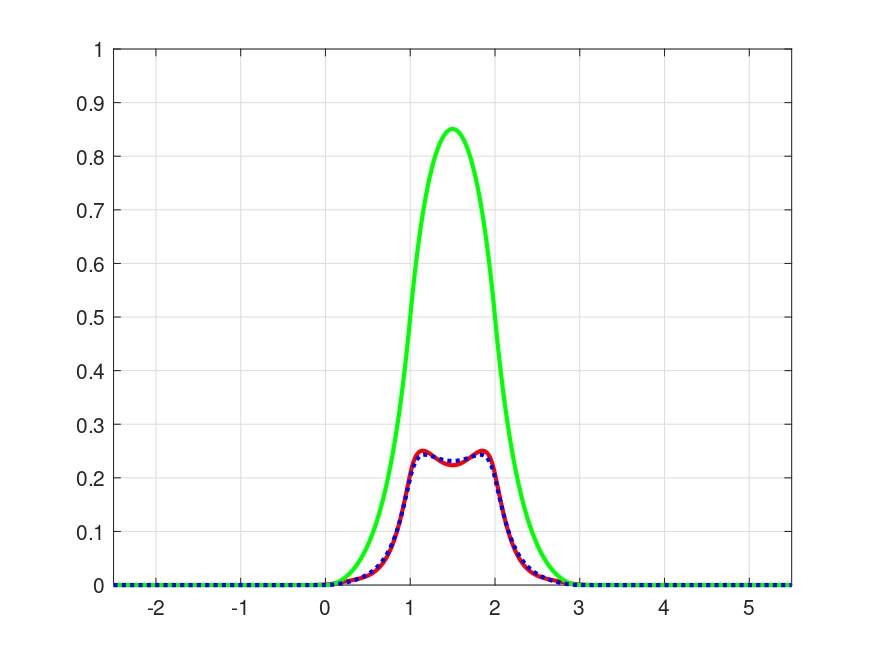}
		\subcaption{}
		\label{fig:sub-c}
	\end{subfigure}
	
	\caption{$a = 1$, $b = \frac{1}{5}$. Green curve represents the generator $\varepsilon_3$. 
		(a) Blue: symmetric dual $k$; red: the dual constructed from $k$ using \ref{stoevaReal}.
		(b) Blue: asymmetric dual $h$; red: the dual constructed from $h$ using \ref{stoevaReal}. 
		(c) Blue: canonical dual $S^{-1}\varepsilon_3$; red: the dual constructed from it using \ref{stoevaReal}.}
	\label{fig:EBS3}
\end{figure}

It can be seen that $
\displaystyle \sum_{n\in\mathbb{Z}} \varepsilon_2'(x-n) = \displaystyle \sum_{n=-1}^{1} \varepsilon_2'(x-n)=\left(\frac{e^{\frac{p}{2}}-e^{-\frac{p}{2}}}{p}\right)^2, \  x \in \mathbb{R}.
$ Then $\varepsilon_2 = \frac{1}{C}\varepsilon_2'$ has partition of unity property, where $C = \left(\frac{e^{\frac{p}{2}}-e^{-\frac{p}{2}}}{p}\right)^2.$
For $\ a_1 = 0, \ a_2 = - a_3= p$ for $p\neq 0,$
\begin{equation*}
	\widetilde{\varepsilon}_3'(x)=
	\begin{cases}
		0, & x \notin \left[-\frac{3}{2},\frac{3}{2}\right], \\[2ex]
		
		\displaystyle 
		\frac{e^{p\left(x+\frac{3}{2}\right)} 
			+ e^{-p\left(x+\frac{3}{2}\right)} -2}{2p^2}, 
		& x \in \left[-\frac{3}{2},-\frac{1}{2}\right], \\[3ex]
		
		\displaystyle 
		\frac{e^{p}+e^{-p}}{p^2}
		-\frac{(1+e^{-p})e^{p\left(x+\frac{1}{2}\right)}}{2p^2}
		-\frac{(1+e^{p})e^{-p\left(x+\frac{1}{2}\right)}}{2p^2},
		& x \in \left[-\frac{1}{2},\frac{1}{2}\right], \\[3ex]
		
		\displaystyle 
		-\frac{1}{p^2}
		+\frac{e^{-p}e^{p\left(x-\frac{1}{2}\right)}}{2p^2}
		+\frac{e^{p}e^{-p\left(x-\frac{1}{2}\right)}}{2p^2},
		& x \in \left[\frac{1}{2},\frac{3}{2}\right].
	\end{cases}
\end{equation*}
Figures \ref{fig:EBS2} and \ref{fig:EBS3} show $\varepsilon_2$ and $\varepsilon_3$ along with their duals obtained from \eqref{stoevaReal}.

It can be verified that the scaled function $ \varepsilon_3 = \frac{ p^2e^p}{(e^p-1)^2}\varepsilon_3' $ satisfies the partition of unity property.\\
In the four cases $B_2$, $B_3$, $\varepsilon_2$, and $\varepsilon_3$, 
we compute the symmetric dual $h$ as defined in \eqref{symdual}, 
the asymmetric dual $k$ given in \eqref{asymdual}, 
and the canonical dual described in \eqref{cano}. 
Using the canonical dual, we further construct alternate duals 
according to  \eqref{altdual*}. 
Additionally, we derive dual windows based on the construction in \eqref{stoevaReal}.  
To construct the dual windows described in \eqref{stoevaReal}, 
for each Gabor window $g = B_2, B_3,$ and $\varepsilon_3$, 
we consider the choices 
$g^d = k$, $h$, $k_2$, $h_2$, and $S^{-1}g$, respectively, 
where 
$
k_2 = S^{-1}g - g + Sk, 
$ and $
h_2 = S^{-1}g - g + Sh.
$
In each case, we fix the Bessel function as 
$
w = \frac{1}{10} g.
$
The resulting dual windows are denoted by 
$\phi_k$, $\phi_h$, $\phi_{k_2}$, $\phi_{h_2}$, 
and $\phi_{S^{-1}g}$, corresponding to the respective choice of $g^d$.

Throughout our experiments, we fix the parameters $a = 1$, 
$b = \frac{1}{5}$, and $p = 3$. 
For the second-order B-splines and exponential B-splines, 
it suffices to consider the integer translations $n = -1, 0, 1$ 
when constructing the dual windows $h$ and $k$, 
as well as when computing $G(x)$. 
For the third-order splines, the translations 
$n = -2, \ldots, 2$ are required. In Figures~\ref{fig:B2}--\ref{fig:EBS3}, we display the windows 
$B_2$, $B_3$, $\varepsilon_2$, and $\varepsilon_3$ (shown in green), respectively. 
Their corresponding dual windows $k$, $h$, and $S^{-1}g$ are plotted in blue. 
The dual windows obtained via \eqref{stoevaReal} are shown in red. 
\begin{figure}[H]
	\centering
	
	\begin{subfigure}[t]{0.32\textwidth}
		\centering
		\includegraphics[width=\textwidth]{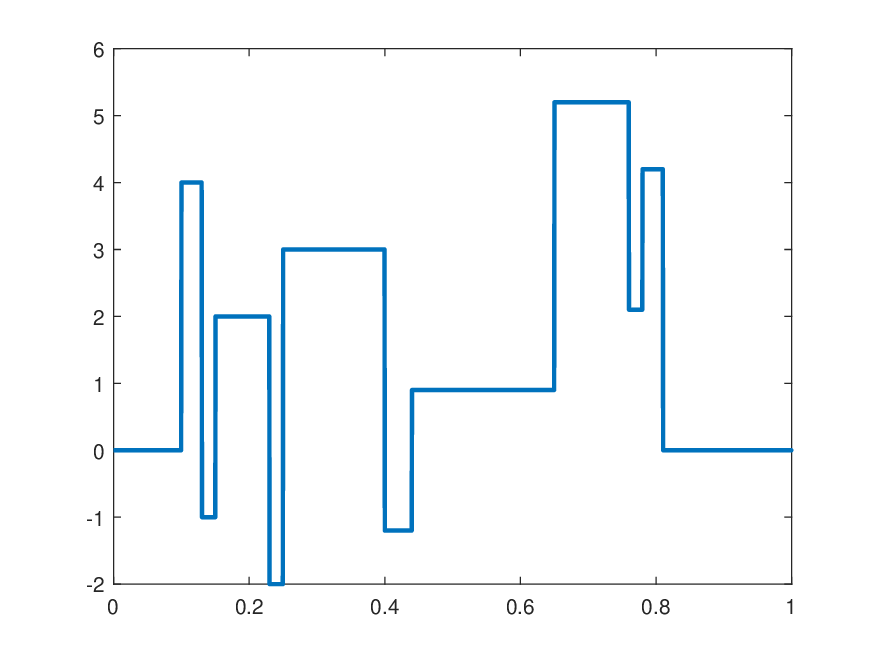}
		\caption{Blocks}
		\label{fig:blocks}
	\end{subfigure}
	\hfill
	\begin{subfigure}[t]{0.32\textwidth}
		\centering
		\includegraphics[width=\textwidth]{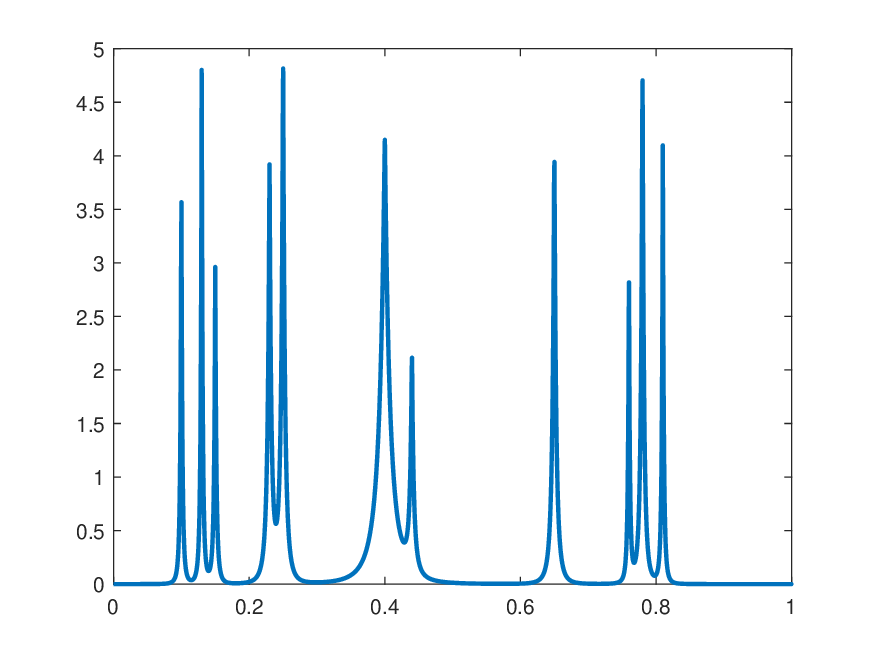}
		\caption{Bumps}
		\label{fig:bumps}
	\end{subfigure}
	\hfill
	\begin{subfigure}[t]{0.32\textwidth}
		\centering
		\includegraphics[width=\textwidth]{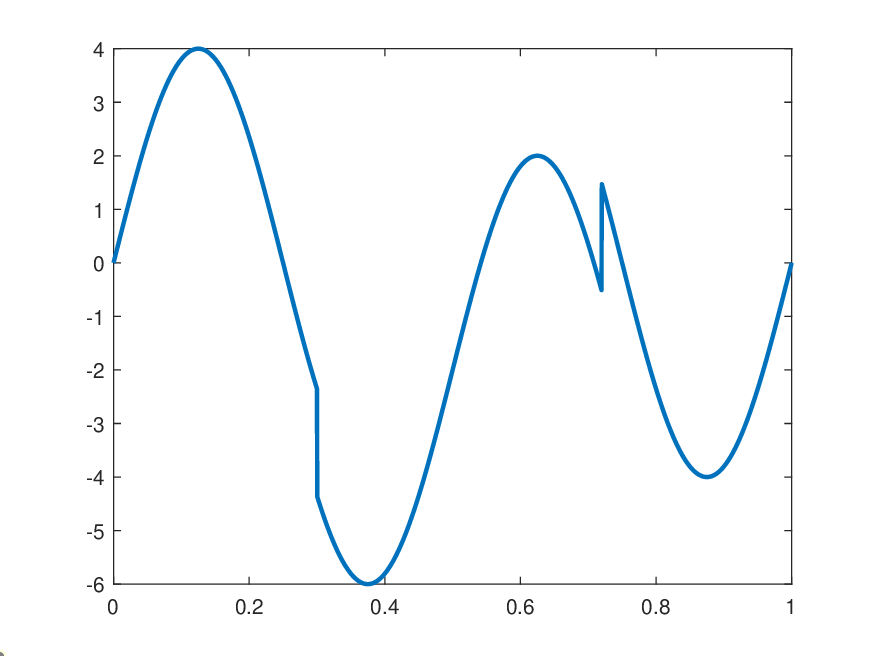}
		\caption{Heavisine}
		\label{fig:heavisine}
	\end{subfigure}
	
	\vspace{1em} 
	
	\begin{subfigure}[t]{0.32\textwidth}
		\centering
		\includegraphics[width=\textwidth]{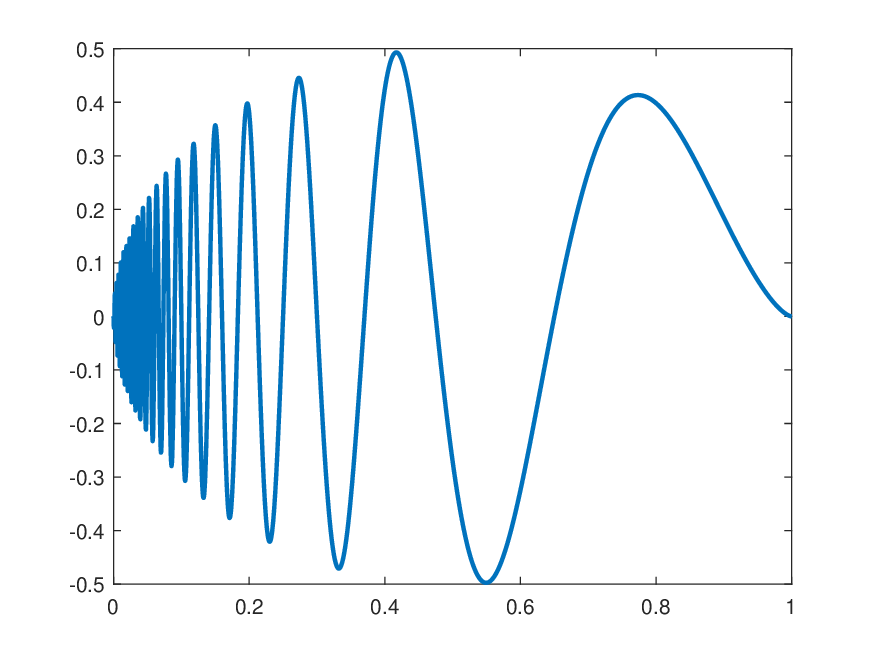}
		\caption{Doppler}
		\label{fig:doppler}
	\end{subfigure}
	\hspace{0.17\textwidth} 
	\begin{subfigure}[t]{0.32\textwidth}
		\centering
		\includegraphics[width=\textwidth]{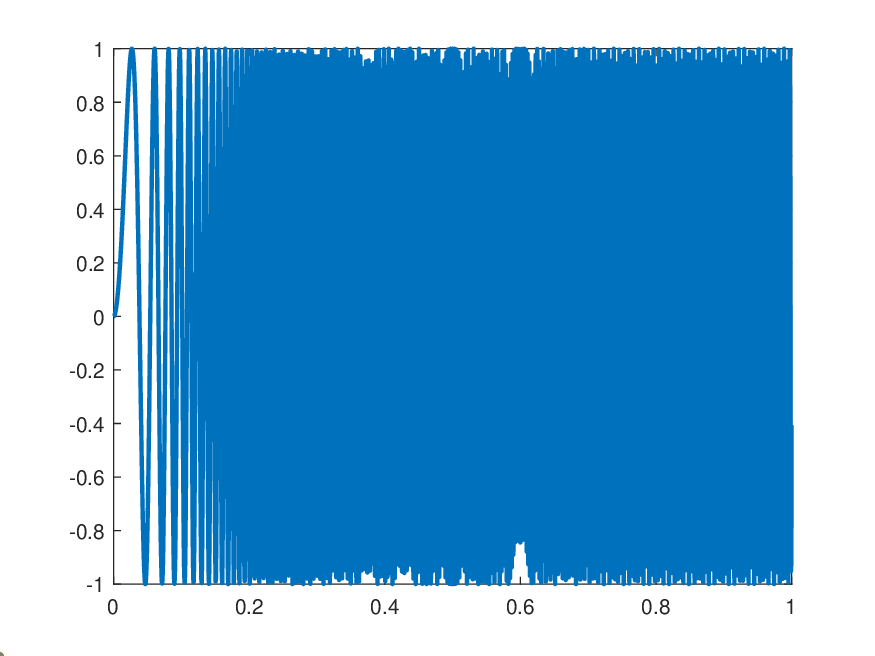}
		\caption{QuadChirp}
		\label{fig:quadchirp}
	\end{subfigure}
	
	\caption{Test signals used in the experiments, each sampled at 2048 points in the interval $[0,1]$.}
	\label{fig:test_signals}
\end{figure}

 \par We use five standard test functions commonly employed in signal processing and denoising Blocks, Bumps, Heavisine, Doppler, and QuadChirp as originally defined in references \cite{don} and \cite{david}. To prepare these functions for numerical analysis, each one is discretized by sampling 2048 equidistant points over the interval 
 $ [0,1] $. This process yields discrete signal vectors that serve as inputs for our experiments.
 Figure 1 shows the resulting signals. 
Each test function $f$ is reconstructed via the Gabor frame expansion
\begin{equation*} \label{frame_exp}
	f = \sum_{m,n \in \mathbb{Z}} \langle f, E_{m b} T_{n a} \mu \rangle\, E_{m b} T_{n a} g 
\end{equation*}
where $\mu$ denotes any of the duals construcetd above. For the numerical computations, we have to fix finite values for $m$ and $n$. We choose $m,n=-3,\dots3.$ The quality of signal reconstruction is quantitatively measured using the Average Mean Square Error (AMSE) metric. The results, summarized in Table~\ref{tab:amse_B2_B3} and Table~\ref{tab:amse_ebs2_ebs3}, show the comparative accuracy of these duals across different signals in \ref{fig:test_signals}  types and generator functions.
 \par
 The following results in the tables show that the duals $ S^{-1}g,\ h,\ k,\  \phi_h,\ \phi_k$ and $\phi_{S^{-1}g}$ achieve almost identical reconstruction across all test signals. 
 In several cases, the symmetric dual $k$ and the dual constructed from it, namely $\phi_k$ yield the smallest AMSE values. On the other hand, the alternative duals $h_2$ and $k_2$ consistently produce larger errors, especially for signals Blocks and Heavisine. Overall, the numerical results demonstrate that the construction \eqref{stoevaReal} preserves, and in some cases slightly improves, reconstruction performance while maintaining stability across different generators and signal classes.
\begin{table}[H]
	\centering
	\renewcommand{\arraystretch}{1.9}
	\resizebox{\textwidth}{!}{%
		\begin{tabular}{|c||ccccc|ccccc|}
			\hline
			\multirow{2}{*}{Dual} & \multicolumn{5}{c|}{$B_2$} & \multicolumn{5}{c|}{$B_3$} \\
			\cline{2-11}
			& Blocks & Bumps & Heav. & Dopp. & Quad 
			& Blocks & Bumps & Heav. & Dopp. & Quad \\
			\hline
			$S^{-1}g$   & 3.4992 & 0.4331 & 8.0253 & 0.0795 & 0.4947 
			& 3.4338 & 0.4313 & 7.9552 & 0.0798 & 0.4947 \\
			h           & 3.3969 & 0.4348 & 7.8954 & 0.0798 & 0.4948 
			& 3.3889 & 0.4345 & 7.8980 & 0.0798 & 0.4948 \\
			k           & 3.3877 & 0.4303 & 7.8626 & 0.0803 & 0.4947 
			& 3.3880 & 0.4303 & 7.8638 & 0.0803 & 0.4947 \\
			$h_2$       & 6.2686 & 0.5329 & 9.7382 & 0.0858 & 0.4959 
			& 6.8218 & 0.5528 & 10.0742 & 0.0861 & 0.4955 \\
			$k_2$       & 3.8635 & 0.4467 & 8.3592 & 0.0867 & 0.4958 
			& 4.2184 & 0.4654 & 8.6151 & 0.0845 & 0.4955 \\
			$\phi_h$    & 3.4246 & 0.4339 & 7.8875 & 0.0797 & 0.4947 
			& 3.3933 & 0.4332 & 7.8971 & 0.0798 & 0.4947 \\
			$\phi_k$    & 3.4225 & 0.4312 & 7.9026 & 0.0800 & 0.4947 
			& 3.3975 & 0.4305 & 7.8910 & 0.0802 & 0.4947 \\
			$\phi_{S^{-1}g}$ & 3.4946 & 0.4329 & 8.0480 & 0.0795 & 0.4947 
			& 3.4319 & 0.4312 & 7.9652 & 0.0798 & 0.4947 \\
			\hline
		\end{tabular}
	}
	\caption{AMSE values for the reconstructed  test signals using Gabor windows $B_2$ and $B_3$ with their respective dual windows.}
	\label{tab:amse_B2_B3}
\end{table}

\begin{table}[H]
	\centering
	\renewcommand{\arraystretch}{1.9}
	\resizebox{\textwidth}{!}{%
		\begin{tabular}{|c||ccccc|ccccc|}
			\hline
			\multirow{2}{*}{Dual} & \multicolumn{5}{c|}{$\varepsilon_2$} & \multicolumn{5}{c|}{$\varepsilon_3$} \\
			\cline{2-11}
			& Blocks & Bumps & Heav. & Dopp. & Quad 
			& Blocks & Bumps & Heav. & Dopp. & Quad \\
			\hline
			$S^{-1}g$   & 3.437 & 0.4297 & 8.0598 & 0.0795 & 0.4947 
			& 3.4852 & 0.43227 & 8.0454 & 0.079622 & 0.4947 \\
			h           & 3.4737 & 0.4345 & 7.873 & 0.0798 & 0.4947 
			& 3.3705 & 0.4367 & 7.9908 & 0.0798 & 0.4948 \\
			k           & 3.3877 & 0.4303 & 7.8626 & 0.0803 & 0.4947 
			& 3.3878 & 0.4303 & 7.863 & 0.0803 & 0.4947 \\
			$h_2$       & 6.2686 & 0.5329 & 9.7382 & 0.0858 & 0.4959 
			& 4.7305 & 0.46243 & 18.6847 & 0.089914 & 0.4951 \\
			$k_2$       & 5.2802 & 0.5228 & 16.8243 & 0.1206 & 0.4951 
			& 4.0635 & 0.4291 & 10.6733 & 0.0857 & 0.4951 \\
			$\phi_h$    & 3.4537 & 0.4324 & 7.9347 & 0.0795 & 0.4947 
			& 3.3906 & 0.4344 & 7.9784 & 0.0798 & 0.4947 \\
			$\phi_k$    & 3.4008 & 0.43025 & 7.9211 & 0.0800 & 0.4947 
			& 3.4069 & 0.4307 & 7.9281 & 0.0800 & 0.4947 \\
			$\phi_{S^{-1}g}$ & 3.6103 & 0.4384 & 8.103 & 0.0793 & 0.4947 
			& 3.479 & 0.4319 & 8.0736 & 0.0797 & 0.4947 \\
			\hline
		\end{tabular}
	}
	\caption{AMSE values for the reconstructed test signals using Gabor windows  $\varepsilon_2$ and $\varepsilon_3$ with their respective dual windows.}
	\label{tab:amse_ebs2_ebs3}
\end{table}

\section{Gabor Frames in Image Reconstruction}

Gabor frames play an important role in the analysis and
reconstruction of images. Due to their time–frequency
localization properties and redundancy, they are particularly
well suited for image reconstruction where stability and robustness
against noise are essential. In order to extend one-dimensional Gabor frame constructions to the
two-dimensional setting required for image processing, we employ the
tensor product structure of Hilbert spaces. A two-dimensional image can be viewed as an element of $L^2(\mathbb{R}^2)$. Since $
L^2(\mathbb{R}^2) \cong L^2(\mathbb{R}) \otimes L^2(\mathbb{R}),
$ two-dimensional Gabor systems can be constructed via tensorization of one-dimensional frames.
If $\mathcal{G}(g_1,a,b)$ and $\mathcal{G}(g_2,a,b)$ are Gabor frames in $L^2(\mathbb{R})$, then we have $\mathcal{G}(g_1 \otimes g_2,a,b):= \{E_{m(b,b)}T_{n(a,a)}g_1\otimes g_2\}
$ is a frame for $L^2(\mathbb{R}^2)$. If $h_1$ and $h_2$ are duals of $g_1$ and $g_2$ respectively, then $h_1 \otimes h_2$ is a dual of $g_1 \otimes g_2$ in $L^2(\mathbb{R}^2)$, see Theorem 2.3 in \cite{arefiimage}. 
This construction provides a direct mechanism for generating
two-dimensional dual Gabor systems from one-dimensional
alternate duals.
 An image $M \in L^2(\mathbb{R}^2)$ (or a discrete
image matrix) can therefore be reconstructed using
\begin{equation}\label{tensor}
M = \sum_{m,n\in\mathbb{Z}} 
\langle M, E_{m(b,b)}T_{n(a,a)}(h \otimes h) \rangle 
E_{m(b,b)}T_{n(a,a)}(g \otimes g),
\end{equation}
where $h$ denotes a chosen dual window. Since $
L^2(\mathbb{R}^2) \cong L^2(\mathbb{R}) \otimes L^2(\mathbb{R})
$, we have $(f \otimes g) (x,y) = f(x)g(y),$ $ \forall x,y \in \mathbb{R}.$ For basic properties of tensor product frames, see \cite{amir}.

In practical computations, only finitely many terms of the sum in \eqref{tensor} are used, leading to an approximation of the original image. The
quality of reconstruction depends strongly on the choice of the dual
window.  Hence, the freedom in selecting dual Gabor windows provides an
additional degree of flexibility in image reconstruction. By
appropriately designing dual generators, one may achieve enhanced
localization, improved numerical stability, and better approximation
performance.

\begin{figure}[htbp]
	\centering
	
	\begin{subfigure}{0.23\textwidth}
		\centering
		\includegraphics[height=3cm]{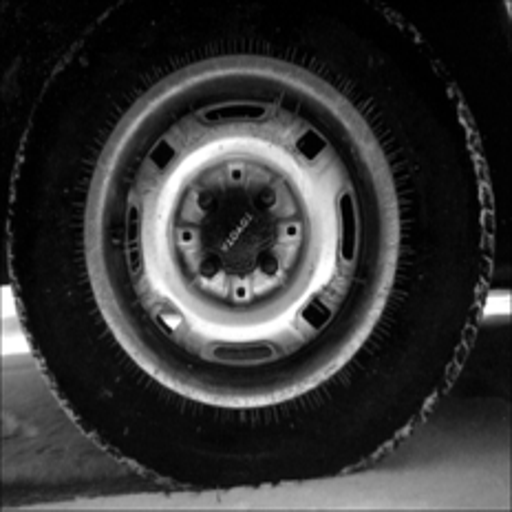}
		\caption{Tire-original}
	\end{subfigure}\hfill
	\begin{subfigure}{0.23\textwidth}
		\centering
		\includegraphics[height=3cm]{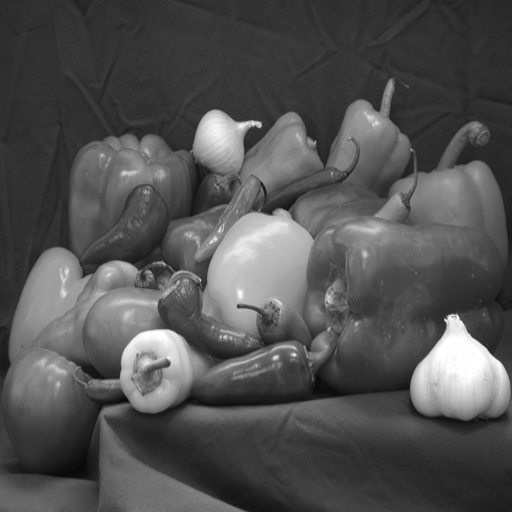}
		\caption{Peppers-original}
	\end{subfigure}\hfill
	\begin{subfigure}{0.23\textwidth}
		\centering
		\includegraphics[height=3cm]{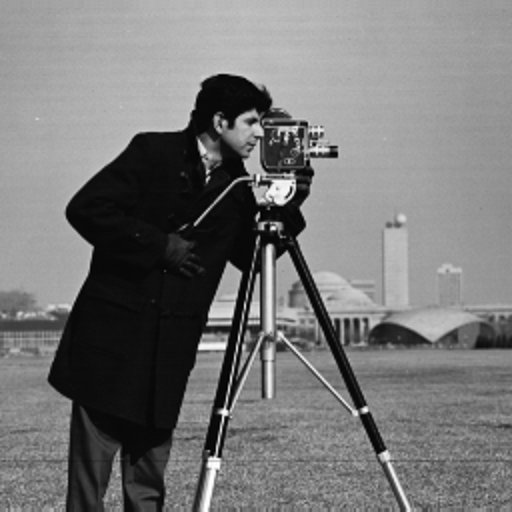}
		\caption{Cameraman-original}
	\end{subfigure}\hfill
	\begin{subfigure}{0.23\textwidth}
		\centering
		\includegraphics[height=3cm]{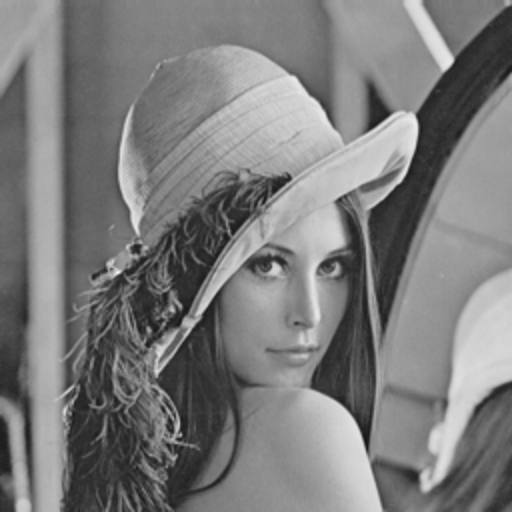}
		\caption{Lena-original}
	\end{subfigure}
	
	\vspace{0.4cm}
	
	\begin{subfigure}{0.23\textwidth}
		\centering
		\includegraphics[width=\linewidth]{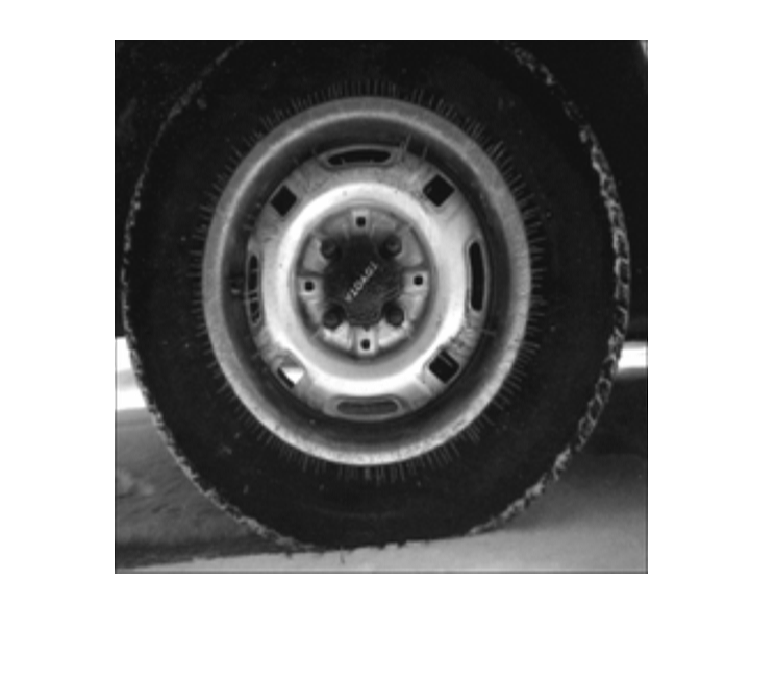}
		\caption{Tire-reconstructed}
	\end{subfigure}\hfill
	\begin{subfigure}{0.23\textwidth}
		\centering
		\includegraphics[width=\linewidth]{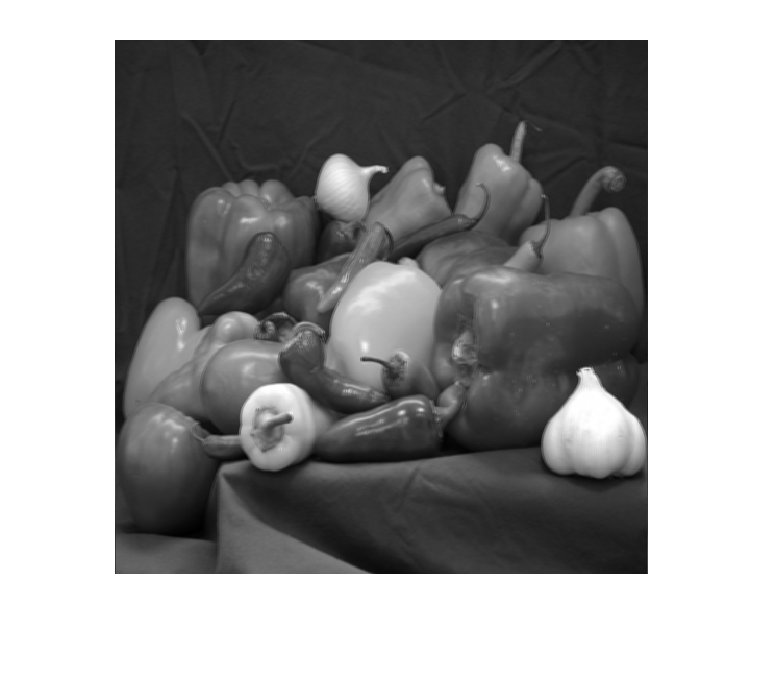}
		\caption{Peppers-reconstructed}
	\end{subfigure}\hfill
	\begin{subfigure}{0.23\textwidth}
		\centering
		\includegraphics[width=\linewidth]{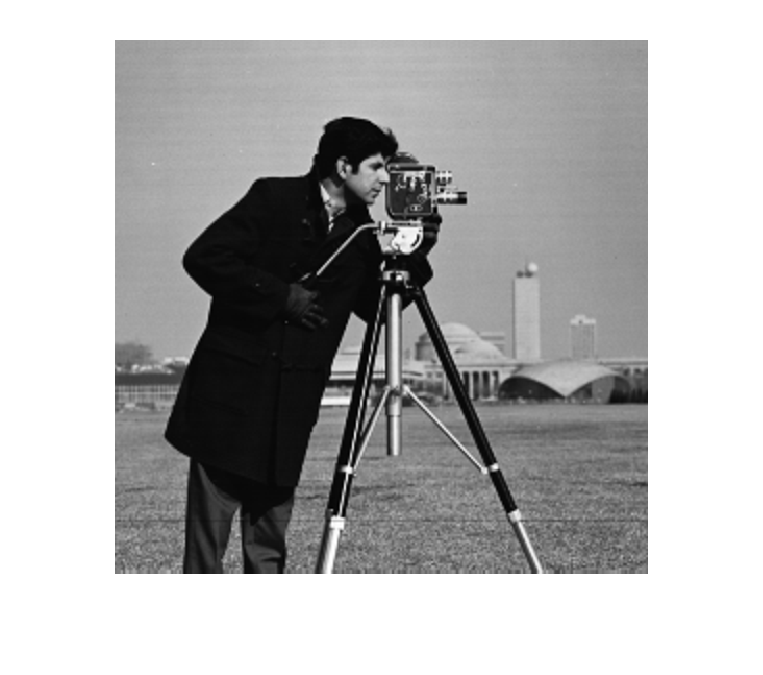}
		\caption{Cameraman-reconstructed}
	\end{subfigure}\hfill
	\begin{subfigure}{0.23\textwidth}
		\centering
		\includegraphics[width=\linewidth]{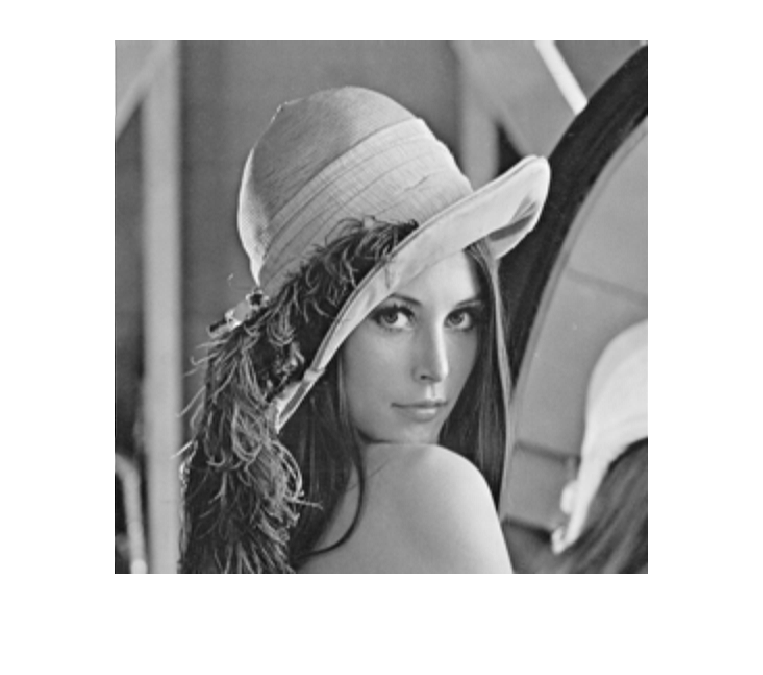}
		\caption{Lena-reconstructed}
	\end{subfigure}
	
	\caption{Reconstruction of test images using different dual windows: First row shows the original images and the second row shows the respective reconstructed images. The image 'Tire' ((a)) is reconstructed with $B_3$ and its dual $\phi_k$, while the image 'Peppers' ((b)) is reconstructed using the dual $k_2$ of $\varepsilon_3$. The image 'Cameraman' ((c)) is reconstructed with $B_2$ and its canonical dual $S^{-1}B_2$, while the image 'Lena' ((d)) is reconstructed using the dual $\phi_{h}$ of $\varepsilon_3$.}
	\label{fig:recon_images}
\end{figure}

 To demonstrate the performance of dual Gabor frames, we reconstruct some test images Cameraman, Lena, Tire and Peppers using the tensor product frames generated by B-splines and exponential B-splines, along with their tensorized duals mentioned above. A finite number of terms in \eqref{tensor} is used to approximate the test image $M$. More precisely, we consider modulation indices $m=-3,\dots 3$ and translation indices 
 $n=-1,\dots W_p$ where $W_p$ denotes the padded image width. The lattice parameters are fixed as 
 $a=1$ and $b= \frac{1}{2N-1}.$ The resulting approximations of some test images obtained from different generator–dual pairs are displayed in Figure~\ref{fig:recon_images}. The performance of each reconstruction is quantified using the average mean-square error (AMSE). Table~\ref{tab:amseBspline} and Table~\ref{tab:amseEBS} summarise the AMSE values corresponding to the different dual Gabor generators. 
 
The AMSE values reported in Tables \ref{tab:amseBspline} and \ref{tab:amseEBS} illustrate the reconstruction accuracy obtained with different dual windows for tensor-product Gabor frames generated by B-splines and exponential B-splines. As expected, the canonical dual $S^{-1}g$ yields extremely small AMSE values, on the order of $10^{-33}$ for $B_2$ and $\varepsilon_2$, which correspond to squared floating-point rounding errors and therefore indicate numerically exact reconstruction in the discrete setting. For the higher-order generators $B_3$ and $\varepsilon_3$, the canonical dual continues to provide the smallest reconstruction error. The alternate compactly supported dual windows produce slightly larger AMSE values, typically around $10^{-5}$, which can be attributed to truncation and finite-support approximation effects rather than instability of the reconstruction. By analysing the AMSE values, it can be observed that the duals $\phi_k$ and $k$ generally yield smaller errors compared to the other compactly supported duals. The alternate duals $h_2$ and $k_2$ tend to produce relatively larger AMSE values. Nevertheless, the reconstruction accuracy remains stable across all test images, indicating that compactly supported dual windows can provide accurate reconstruction.
\begin{table}[htbp]
	\centering
	\renewcommand{\arraystretch}{1.4}
	\resizebox{\textwidth}{!}{%
		\begin{tabular}{|c||cccc|cccc|}
			\hline
			\multirow{2}{*}{Dual} 
			& \multicolumn{4}{c|}{$B_2$} 
			& \multicolumn{4}{c|}{$B_3$} \\
			\cline{2-9}
			& Tire & Peppers & Cameraman & Lena 
			& Tire & Peppers & Cameraman & Lena \\
			\hline
			
			$S^{-1}g$   
			& $1.59\times10^{-32}$ & 
			$2.45\times10^{-33}$ & $1.81\times10^{-30}$ & $8.66\times10^{-31}$
			& $2.41\times10^{-6}$ & $4.03\times10^{-6}$ & $9.34\times10^{-6}$ & $5.51\times10^{-6}$ \\
			
			$h$      
			& $4.93\times10^{-4}$ & $4.18\times10^{-4}$ & $1.40\times10^{-3}$ & $8.38\times10^{-4}$
			& $2.91\times10^{-5}$ & $5.09\times10^{-5}$ & $1.13\times10^{-4}$ & $6.71\times10^{-5}$ \\
			
			$k$         
			& $6.91\times10^{-5}$ & $1.15\times10^{-4}$ & $2.79\times10^{-4}$ & $1.59\times10^{-4}$
			& $1.52\times10^{-5}$ & $2.17\times10^{-5}$ & $2.79\times10^{-5}$ & $3.19\times10^{-5}$ \\
			
			$h_2$       
			& $8.68\times10^{-4}$ & $7.81\times10^{-4}$ & $2.60\times10^{-3}$ & $1.50\times10^{-3}$
			& $6.25\times10^{-4}$ & $1.00\times10^{-3}$ & $2.70\times10^{-3}$ & $1.40\times10^{-3}$ \\
			
			$k_2$       
			& $1.49\times10^{-4}$ & $2.50\times10^{-4}$ & $5.71\times10^{-4}$ & $3.41\times10^{-4}$
			& $3.05\times10^{-4}$ & $3.54\times10^{-4}$ & $7.58\times10^{-4}$ & $6.26\times10^{-4}$ \\
			
			$\phi_h$    
			& $8.66\times10^{-5}$ & $7.22\times10^{-5}$ & $2.51\times10^{-4}$ & $1.46\times10^{-4}$
			& $1.43\times10^{-5}$ & $2.55\times10^{-5}$ & $5.18\times10^{-5}$ & $3.28\times10^{-5}$ \\
			
			$\phi_k$    
			& $1.21\times10^{-5}$ & $1.99\times10^{-5}$ & $5.01\times10^{-5}$ & $2.78\times10^{-5}$
			& $8.15\times10^{-6}$ & $1.27\times10^{-5}$ & $1.42\times10^{-5}$ & $1.74\times10^{-5}$ \\
			
			$\phi_{S^{-1}g}$ 
			& $4.94\times10^{-8}$ & $8.17\times10^{-8}$ & $2.00\times10^{-7}$ & $1.13\times10^{-7}$
			& $2.36\times10^{-6}$ & $3.96\times10^{-6}$ & $9.18\times10^{-6}$ & $5.42\times10^{-6}$ \\
			
			\hline
		\end{tabular}%
	}
	\caption{AMSE values for the test images reconstructed using tensor-product Gabor frames generated by $B_2$ and $B_3$ with their respective dual windows.}
	\label{tab:amseBspline}
\end{table}
\begin{table}[htbp]
	\centering
	\renewcommand{\arraystretch}{1.6}
	\resizebox{\textwidth}{!}{%
		\begin{tabular}{|c||cccc|cccc|}
			\hline
			\multirow{2}{*}{Dual} 
			& \multicolumn{4}{c|}{$\varepsilon_2$} 
			& \multicolumn{4}{c|}{$\varepsilon_3$} \\
			\cline{2-9}
			& Tire & Peppers & Cameraman & Lena 
			& Tire & Peppers & Cameraman & Lena \\
			\hline
			
			$S^{-1}g$   
			& $1.59\times10^{-32}$ & 
		$2.45\times10^{-33}$ & $1.81\times10^{-30}$ & $8.66\times10^{-31}$
			& $6.62\times10^{-7}$ & $1.10\times10^{-6}$ & $2.62\times10^{-6}$ & $1.52\times10^{-5}$ \\
			
			$h$      
			& $4.93\times10^{-4}$ & $4.18\times10^{-4}$ & $1.40\times10^{-3}$ & $8.38\times10^{-4}$
			& $3.55\times10^{-5}$ & $5.86\times10^{-5}$ & $1.41\times10^{-4}$ & $8.06\times10^{-5}$ \\
			
			$k$         
			& $6.91\times10^{-5}$ & $1.15\times10^{-4}$ & $2.79\times10^{-4}$ & $1.59\times10^{-4}$
			& $1.52\times10^{-5}$ & $2.12\times10^{-5}$ & $2.66\times10^{-5}$ & $3.18\times10^{-5}$ \\
			
			$h_2$       
			& $1.20\times10^{-3}$ & $1.10\times10^{-3}$ & $3.50\times10^{-3}$ & $2.00\times10^{-3}$  
			& $9.30\times10^{-3}$ & $1.46\times10^{-2}$ & $3.74\times10^{-2}$ & $2.12\times10^{-2}$ \\
			
			$k_2$       
			& $1.12\times10^{-4}$ & $1.88\times10^{-4}$ & $4.35\times10^{-4}$ & $2.57\times10^{-4}$  
			& $4.20\times10^{-3}$ & $5.10\times10^{-3}$ & $9.00\times10^{-3}$ & $8.60\times10^{-3}$ \\
			
			$\phi_h$    
			& $1.10\times10^{-4}$ & $9.11\times10^{-5}$ & $3.19\times10^{-4}$ & $1.85\times10^{-4}$ 
			& $1.24\times10^{-5}$ & $2.08\times10^{-5}$ & $4.70\times10^{-5}$ & $2.81\times10^{-5}$ \\
			
			$\phi_k$    
			& $1.25\times10^{-5}$ & $2.05\times10^{-5}$ & $5.14\times10^{-5}$ & $2.86\times10^{-5}$  
			& $5.68\times10^{-6}$ & $8.53\times10^{-6}$ & $9.38\times10^{-6}$ & $1.20\times10^{-5}$ \\
			
			$\phi_{S^{-1}g}$ 
			& $3.47\times10^{-7}$ & $2.79\times10^{-7}$ & $9.93\times10^{-7}$ & $5.76\times10^{-7}$  
			& $5.96\times10^{-7}$ & $9.92\times10^{-7}$ & $2.37\times10^{-6}$ & $1.37\times10^{-6}$ \\
			
			\hline
		\end{tabular}
	}
	\caption{AMSE values for the test images reconstructed using tensor-product Gabor frames generated by $\varepsilon_2$ and $\varepsilon_3$ with their respective dual windows.}
	\label{tab:amseEBS}
\end{table}
\section{Conclusion}
In this work, we investigated the performance of several compactly supported dual window constructions for Gabor frames generated by B-splines and exponential B-splines. Starting from explicit constructions of symmetric and asymmetric duals, we further examined alternate duals obtained via iterative perturbation methods and through the general characterization of compactly supported dual frames.  The numerical experiments for both one-dimensional test signals and two-dimensional image reconstruction reveal several important observations. First, symmetric compactly supported duals frequently achieve reconstruction errors comparable to the canonical dual, despite avoiding operator inversion. Second, duals constructed using \eqref{stoevaReal} preserve stability while maintaining compact support, making them computationally attractive. Third, exponential B-spline generators consistently yield lower AMSE values across a wide range of signals and images, indicating improved reconstruction capability relative to classical polynomial B-splines.

In short, we observe that compactly supported dual windows are not only theoretically viable but also practically competitive. The results suggest that symmetric compactly supported duals, particularly when combined with exponential B-spline generators, provide an effective balance between localization, computational efficiency, and reconstruction accuracy. These findings contribute to the practical design of Gabor frame systems in signal and image processing applications.

	\end{document}